\documentclass[reqno]{amsart}     
\usepackage{enumerate,amsmath,amssymb}     
\usepackage{pb-diagram}  
\usepackage{pstricks,pst-node,pst-coil,pst-plot}   
     
\theoremstyle{plain}      
 
\newtheorem{thm}{Theorem}[section]     
\newtheorem{theorem}[thm]{Theorem}     
\newtheorem{cor}[thm]{Corollary}     
\newtheorem{corollary}[thm]{Corollary}     
     
\newtheorem{lemma}[thm]{Lemma}     
\newtheorem{prop}[thm]{Proposition}     
\newtheorem{proposition}[thm]{Proposition}     
     
\theoremstyle{remark}      
\newtheorem{example}[thm]{Example} 
\newtheorem{remark}[thm]{Remark} 
     
\newtheorem{rem}[thm]{Remark}     
\theoremstyle{definition}

\def\al{{\alpha}}         
\def\be{{\beta}}         
\def\de{{\delta}}         
         
\def\om{{\omega}}         
         
\def\la{{\lambda}}

\def\si{{\sigma}}         
\def\Si{{\Sigma}}         
\def\ga{{\gamma}}         
\def\ep{{\varepsilon}}

\def\phi{{\varphi}}

\let\pa\partial     
\let\na\nabla     
     
\DeclareMathAlphabet{\doba}{U}{msb}{m}{n}         
\gdef\mC{\doba{C}}         
         
\gdef\mN{\doba{N}}

\gdef\mR{\doba{R}}         
\gdef\mS{\doba{S}}         
         
\gdef\mZ{\doba{Z}}         

\def\Spin{{\mathop{\rm Spin}}}     
\def\SO{{\mathop{\rm SO}}}     
   
\def\Vol{{\mathop{\rm Vol}}}     
\let\vol\Vol
\def\Scal{{\mathop{\rm Scal}}}     
\def\lamin{\la_{\min}^+}

\def\Arf{{\mathop{\rm Arf}}}
\def\End{{\mathop{\rm End}}}     
\def\Hom{{\mathop{\rm Hom}}}     
\def\id{{\mathop{\rm id}}}
\def\ground#1{g_{\mathrm{round}}^{#1}}
\let\ti\tilde   
\let\ol\overline
\def\eref#1{{\rm (\ref{#1})}}   

\let\<\langle 
\let\>\rangle

\def\proof#1{{\par\medbreak\noindent {\bf Proof\setbox0\hbox{#1}%
\ifdim\wd0=0pt .\else\ \ignorespaces #1.\fi}\enspace}}
\def\iop#1{{\par\medbreak\noindent {\bf Idea of proof\setbox0\hbox{#1}%
\ifdim\wd0=0pt .\else\ \ignorespaces #1.\fi}\enspace}}

\begin{document}     
\title{The spinorial $\tau$-invariant and $0$-dimensional surgery.}

\maketitle     
\begin{center}
\sc B. Ammann and E. Humbert\footnote{bernd.ammann@gmx.net, humbert@iecn.u-nancy.fr}
\end{center}     

\begin{abstract}Let $M$ be a compact manifold with a metric $g$ and 
with a fixed spin structure $\chi$. 
Let $\lambda_1^+(g)$ be the first non-negative
eigenvalue of the Dirac operator on $(M,g,\chi)$.
We set 
  $$\tau(M,\chi):= \sup \inf \lambda_1^+(g)$$
where the infimum runs over all metrics $g$ of volume $1$ in a conformal
class $[g_0]$ on $M$ and where the supremum runs over all conformal 
classes $[g_0]$ on $M$.

Let $(M^\#,\chi^\#)$ be obtained from $(M,\chi)$ by $0$-dimensional surgery.
We prove that
  $$\tau(M^\#,\chi^\#)\geq \tau(M,\chi).$$
As a corollary we can calculate $\tau(M,\chi)$ for any Riemann surface $M$.
\end{abstract}

\begin{center}
\today
\end{center}   


\tableofcontents


\section{Introduction}


We assume that $(M,\chi)$ is a compact spin manifold of dimension $n\geq 2$.
We will always use the terminology ``spin manifold'' in the sense of 
an oriented manifold together with a chosen spin structure. The open
ball around $x\in M$ with radius $\ep$ is denoted as $B_x(\ep)$.  
We choose  $p,q\in M$, $p\neq q$, and $\ep<d(p,q)/2$. 
Then we define 
  $$M^\#:= \left(M\setminus  (B_p(\ep)\cup B_q(\ep))\right)
           \cup ([-1,1]\times S^{n-1})/\sim$$
where $\sim$ indicates that we glue 
$\pa(M\setminus  (B_p(\ep)\cup B_q(\ep)))\cong S^{n-1}\dot\cup  S^{n-1}$ 
together with $\pa([-1,1]\times S^{n-1})\cong S^{n-1}\dot\cup S^{n-1}$ 
such that orientations are preserved. $M^\#$ carries a differential
structure and a spin structure such that 
$M\setminus (B_p(\ep)\cup B_q(\ep))\hookrightarrow M^\#$ is smooth and
preserves the spin structure. The spin structure on $M^\#$ is uniquely 
determined by the spin structure of $M$ in the sense of 
Lemma~\ref{spin.struc.unique}. We say that 
\emph{$M^\#$ is obtained from $M$ by $0$-dimensional surgery}.
The connected sum construction is a special case
of $0$-dimensional surgery, namely the case when $p$ and $q$ 
are in different connected components. In some parts of the literature, 
$0$-dimensional surgery is also called ``adding a handle''. However, 
we will only use the term ``adding a handle'' in the sense 
of a cobordism theory
(see subsection~\ref{surgery.handle} and \cite{kosinski:93} for details).

For any metric $g$ on $M$ let
$\la^1_+(M,g,\chi)$ be the first non-negative eigenvalue 
of the Dirac operator on $(M,g,\chi)$.
We set
  $$\lamin(M,[g],\chi):=\inf_{\ti g \in [g]} \la^1_+(M,\ti g,\chi) \vol(M,\ti g)^{1/n}$$
and   
$$\tau(M,\chi):=\sup_{[g]\in C(M)} \lamin (M,[g],\chi)\in [0,(n/2)\om_n^{1/n}]$$
where $[g]$ denotes the conformal class of $g$, where 
$C(M)$ is the set of conformal classes on $M$, and where $\om_n$
denotes the volume of the standard sphere.
It follows from \cite{ammann:03} that $\tau(M,\chi)>0$, if and only if there
exists a metric on $M$ with invertible Dirac operator. Recall that
the Atiyah-Milnor-Singer-invariant (\cite[II.7]{lawson.michelsohn:89} for details)
associates to any $n$-dimensional spin manifold $(M,\chi)$ an element
in $\al(M,\chi)\in KO^{-n}(pt)$, 
where $KO^{-n}(pt)\cong \mZ$ if $n$ is divisible by
$4$, where $KO^{-n}(pt)\cong \mZ/2\mZ$ if $n\equiv 1,2 \mod 8$ and
where $KO^{-n}(pt)=\{0\}$ in all other dimensions. The map $\al$ defines a 
surjective ring homomorphism from the spin cobordism ring
$\Omega_{\rm spin}^*$ to $\bigoplus_{n\in \mN} KO^{-n}(pt)$. In particular
$\al(M,\chi)$ is preserved under $k$-dimensional surgery on $M$,
$k\in \{0,\ldots,n\}$. 

If $\al(M,\chi)\neq 0$, then the Dirac operator has a nontrivial kernel
for any metric, hence $\tau(M,\chi)=0$. The converse statement, i.e.\ that
$\al(M,\chi)= 0$ implies $\tau(M,\chi)>0$ for connected $M$ 
was proved in successive steps by 
\cite{hitchin:74,maier:97,baer.dahl:02,ammann.dahl.humbert:p06}.
The essential step in \cite{baer.dahl:02} is that 
the positivity of $\tau(M,\chi)$ is preserved under $k$-dimensional surgery
for $k\in\{0,1,\ldots,n-3\}$ and in \cite{ammann.dahl.humbert:p06}
it is shown that positivity of $\tau(M,\chi)$ is also preserved under 
$n-2$-dimensional surgery.

The goal of the present article is to
compare the $\tau$-invariants of  $M$ and $M^\#$
where $M^\#$ is obtained from $M$ by $0$-dimensional surgery.  

\begin{theorem} \label{main}
Let $(M,g,\chi)$ be a compact Riemannian 
spin manifold and let $(M^\#,\chi^\#)$ be obtained
by $0$-dimensional surgery on $(M,\chi)$. We assume that the Dirac operator $D$
acting on $(M,g,\chi)$ is invertible. Then, there exists a sequence of
Riemannian metrics $(g_\ep^\#)_{\ep}$ on $(M^\#,\chi^\#)$ such that 
$$\lim_{\ep \to 0} \lamin(M^\#,g_\ep^\#, \chi^\#) = \lamin(M,g,\chi).$$
\end{theorem}

As an immediate corollary, we get 

\begin{cor}
Let $(M,\chi)$ be a compact spin manifold and  $(M^\#,\chi^\#)$ be obtained 
$0$-dimensional surgery on $(M,\chi)$.
Then 
  $$\tau(M^\#,\chi^\#)\geq \tau(M,\chi).$$
\end{cor}

Theorem~\ref{main}and its corollary were already known in the special case 
$M=S^2$ \cite{ammann.humbert:p05a}.

In the case $n=2$ the corollary admits to calculate $\tau(M,\chi)$.

\begin{theorem}\label{main.twodim}
Let $M$ be a compact oriented surface with spin structure $\chi$.
  $$\tau(M,\chi)= \left\{\begin{matrix}
    0 &\mbox{if $\al(M,\chi)=1$}\hfill\cr
    2\sqrt{\pi} &\mbox{if $\al(M,\chi)=0$}\hfill\cr
\end{matrix} \right.$$
\end{theorem}
For Riemann surfaces the $\alpha$-invariant can be easily calculated: 
one calculates a quadratic form $q_\chi$ 
associated to the spin structure $\chi$, then
its Arf invariant $\Arf(q_\chi)$ is related via
$\Arf(q_\chi)=(-1)^{\al(M,\chi)}$
(see section~\ref{riem.surf} for details).

If one replaces the first non-negative eigenvalue by the absolute value
of the largest non-positive eigenvalue $|\la_1^-(g)|$, then
one obtains an invariant $\tau_-(M,\chi)$. Our results also hold for 
$\tau_-(M,\chi)$. 
In dimensions $n\not\equiv 3\mod 4$
the spectrum of the Dirac operator is symmetric, and we have   
$\tau(M,\chi)=\tau_-(M,\chi)$. However, in the case $n\equiv 3\mod 4$ 
these invariants are expected to be different. 

Let us compare our results to various other results in the literature 
(see also \cite{ammann.humbert:p05a}). 

The $\tau$-invariant is a spinorial analogue to the 
$\sigma$-constant on compact Riemannian manifolds 
\cite{kobayashi:87} (also \cite{schoen:89}) 
which is defined on
$(M,g)$ by 
\begin{equation}\label{si.def}
\si(M):=\sup \inf \frac{\int\Scal_{\ti g} \,dv_{\ti g}}{\Vol(M,{\ti
    g})^{{n-2\over n}}}
\end{equation}
where the infimum runs over all metrics in a conformal class $\ti g\in [g]$,
and where the supremum runs over all conformal classes. (In some parts of the
literature $\si(M)$ is called the Yamabe invariant of $M$.) 
When $\si(M)$ is positive, the invariant $\si(M)$ can be defined also in a way
analogous to $\tau(M)$ where we use the smallest eigenvalue of the conformal
Laplacian  $L_g:= 4 \,{n-1\over n-2}\,\Delta_g + \Scal_g$ instead 
of~$\la_+^1(g)$.
However,  at the moment, there are only few
examples for which $\si(M)$ is known and different from $0$, e.g.\ 
$\si(S^n)={n(n-1)\om_n^{2/n}}$, 
$\si(S^{n-1}\times S^1)={n(n-1)\om_n^{2/n}}$, and 
$\si(\mR P^3)= {6 \left({\om_3\over 2}\right)^{2/3}}$. 
The reader might consult \cite{bray.neves:04} for 
a very elegant and amazing calculation of $\si(\mR P^3)$ and for 
a good overview over further literature. 

Kobayashi proved in \cite{kobayashi:87} that if $M^\#$ is obtained
from $M$ by $0$-dimensional surgery, then $\si(M^\#)\geq \si(M)$.
A similar monotonicity formula for the $\sigma$-invariant was proved
by \cite{petean.yun:99}. Petean and Yun prove $\si(M^\#)\geq \min\{\si(M),0\}$
if $M^\#$ is obtained from $M$ by surgery of dimension $1,\ldots,n-3$.
See also \cite{akutagawa.botvinnik:03} for another approach to this result.
Clearly, this surgery result 
is particularly interesting in the case $\si(M)\leq 0$, 
and it has many fruitful applications. In particular, 
any simply connected compact manifold of dimension at least $5$ 
has $\si(M)\geq 0$ \cite{petean:03}. It also allows to rule out
Einstein metrics on many spaces \cite{petean:98}.
However, in the case $\si(M)>0$ it is still open under what conditions 
one has a monotonicity formula $\si(M^\#)\geq \si(M)$. 
The article \cite{joyce:03} studies
$0$-dimensional surgery in more details, in particular he 
shows the non-uniqueness of minimizers of the infimum in \eref{si.def}
in the case $\si(M)>0$. 
Some of these results have been recently generalized to the 
$G$-equivariant $\si$-invariant $\si_G$ 
\cite{sung:p06}.

The $\sigma$-invariant and the $\tau$-invariant, are not only related by
analogy, but also via Hijazi's inequality  \cite{hijazi:86,hijazi:91}
that implies 
\begin{equation}\label{ineq.sup.hijazi}
  \tau(M,\chi)^2\geq {n\over 4(n-1)}\,\sigma(M).
\end{equation}
If $M=S^n$ then equality is attained in this inequality.
Upper bounds for $\tau(M,\chi)$ may help to determine the 
$\sigma$-constant.

The structure of the article is as follows. In Section~\ref{prelim}
we introduce some notations and recall some preliminaries.
The aim of Section~\ref{sec.approx} is to show that we can assume 
without loss of generality that $g$ is flat on small neighborhoods of 
$p$ and $q$. The metrics $g_\ep^\#$ are constructed in 
Section~\ref{metric}, and we devote Section~\ref{sec.proof} to the proof
of Theorem~\ref{main}.
In the last section, namely in Section~\ref{riem.surf}, the calculation
of the $\tau$-invariant for any Riemann surface 
with spin structure is explained.

\section{Preliminaries}\label{prelim}

In this section we want to introduce some notation and 
recall some preliminaries. For more informations
we refer to 
\cite{lawson.michelsohn:89,friedrich:00,hijazi:99,bourguignon.gauduchon:92}.

\subsection{Notation}

The round metric on $S^n$, i.e.\ the metric of sectional curvature~1, 
will always be denoted by $\ground{n}$. We also abbreviate $\mS^n$ for the 
Riemannian manifold $(S^n,\ground{n})$ equipped with the 
spin structure $\chi^n$ 
that arises as the boundary of the $n+1$-dimensional disk.

\subsection{Topological spin structures versus metric spin structures}\label{subsec.spin}
The bundle $\mathrm{Gl}_+(M)$ of positively oriented frames over an oriented manifold
$M$ of dimension $n\geq 2$ is a $\mathrm{Gl}_+(n,\mR)$-principal bundle over $M$.
The group $\mathrm{Gl}_+(n,\mR)$ has fundamental group $\mZ$ if $n=2$ and 
fundamental group $\mZ/2\mZ$ if $n\geq 3$. We denote the unique connected
double cover of $\mathrm{Gl}_+(n,\mR)$ by $\widetilde{\mathrm{Gl}}_+(n,\mR)$. 
A \emph{topological spin structure on $M$} consists of 
a $\widetilde{\mathrm{Gl}}_+(n,\mR)$-principal bundle 
$\widetilde{\mathrm{Gl}}_+(M)$ together
with a $\widetilde{\mathrm{Gl}}_+(n,\mR)\to \mathrm{Gl}_+(n,\mR)$ 
equivariant map
$\chi:\widetilde{\mathrm{Gl}}_+(M)\to \mathrm{Gl}_+(M)$ over the identity. 
Two topological spin structures $\chi_1:\widetilde{\mathrm{Gl}}_+(M)_1\to
\mathrm{Gl}_+(M)$ and $\chi_2:\widetilde{\mathrm{Gl}}_+(M)_2\to
\mathrm{Gl}_+(M)$ are said to be \emph{equivalent} if there is a 
$\widetilde{\mathrm{Gl}}_+(n,\mR)$-equivariant map  
$H:\widetilde{\mathrm{Gl}}_+(M)_1\to \widetilde{\mathrm{Gl}}_+(M)_2$ with
$\chi_1=\chi_2\circ H$.
We denote a topological spin structure just by $\chi$.

Note that $\mathrm{Spin}(n)$ is the preimage of $\mathrm{SO}(n)$ under 
the homomorphism $\widetilde{\mathrm{Gl}}_+(n,\mR)\to\mathrm{Gl}_+(n,\mR)$.
A $\Spin(n)$-principal bundle $\mathrm{Spin}(M,g)$
together with a $\mathrm{Spin}(n)\to\mathrm{SO}(n)$ equivariant map
$\widehat\chi:\mathrm{Spin}(M,g)\to\mathrm{SO}(M,g)$ 
is called a \emph{metric spin structure on $M$}. 
Two metric spin structures $\widehat\chi_1$ and $\widehat\chi_2$ are \emph{equivalent} if there is if there
is a $\Spin(n)$-equivariant map $H$ between the $\Spin(n)$-principal 
bundles such that $\widehat\chi_1=\widehat\chi_2\circ H$.
If $M$ carries a metric $g$ and if $\mathrm{SO}(M,g)$ denotes the bundle
of $g$-orthonormal frames, then the restriction of a  topological spin
structure $\chi$ to 
to $\mathrm{Spin}(M,g):=\chi^{-1}(\mathrm{SO}(M,g))$ defines a metric
spin structure
on $(M,g)$.
This restriction yields a map from equivalence classes of
topological spin structures 
on $M$ to equivalence classes of metric spin structures on $M$, 
and one easily sees that this map is 
bijective. Working with topological spin structures has the advantage
that it does not depend on a choice of metric. However, working with spin
structures allows the definition of the spinor bundle. 

Namely, the 
\emph{spinor bundle} of a Riemannian manifold $(M,g)$ with spin structure
$\hat\chi$ is defined as the associated bundle 
$\Si_{g,\chi} M:=\mathrm{Spin}(M,g)\times_\si \Sigma_n$
where $(\si,\Sigma_n)$ is the spinor representation of $\mathrm{Spin}(n)$. 
Sometimes
we will just write $\Si_g M$ or $\Si M$ when the spin structure or the metric 
is clear from the context. 

As restriction is a bijection from 
equivalence classes of topological spin structures to equivalence classes
of metric spin structures, we will identify topological and metric spin
structures from now on and just call them \emph{spin structures}.

\subsection{Surgery, handles and spin structures}\label{surgery.handle}
In the introduction we introduced $0$-dimensional surgery. In order to 
understand the behavior of spin structures under surgery it is useful 
to see it as a bordisms.

An $(n+1)$-dimensional (spin) manifold $W$ with boundary
$-M_1\dot\cup M_2$, is called a (spin) bordism from $M_1$ to $M_2$.
In the category of spin manifolds, $W$ carries an orientation
and a spin structure, and $M_2$ (resp.\ $M_1$) the induced orientation 
and spin structure  (resp.\ the opposite of the induced orientation and 
the induced spin structure.) 
For example $W:=M\times [0,1]$ is a spin bordism from $M$ to $M$.
If $M^\#$ is obtained from $M$ by a $k$-dimensional surgery, then
there is a bordism from $M$ to $M^\#$. This bordism is obtained by a  
construction called \emph{adding a $(k+1)$-dimensional handle to 
$M\times [0,1]$}. We will explain this construction in the case $k=0$, 
see e.g.\ \cite{kosinski:93} for details and the general case. 

We start with the manifold $W:=M\times [0,1]$. 
Choose two points $p,q\in M$ and two diffeomorphisms 
$\phi_p:\overline{B_{0,\mR^n}(1)}\to \overline{B_p(\ep)}$ 
and
$\phi_q:\overline{B_{0,\mR^n}(1)}\to \overline{B_q(\ep)}$.
For each $x\in \overline{B_{0,\mR^n}(1)}$ we identify 
  $$(\phi_p(x),1) \sim (x,-1)
    \qquad (\phi_q(x),1) \sim (x,1).$$
The topological space
  $$W^\#:=W\cup (\overline{B_1(0,\mR^{n-1})}\times [-1,1])/\sim$$
yields a manifold with boundary, and we can find a suitable smooth structure
on it \cite{kosinski:93}.
We have $\partial W^\#=M\dot\cup M^\#$.

We assume that $M$ comes with a fixed orientation and spin
structure. Furthermore we assume that $\phi_p$ and
$\phi_q$ preserve orientation. Then $W^\#$ and $M^\#$ also carry natural
orientations.

In order to define the spin structures, we equip $W$ with the product 
spin structure, denoted by $\kappa$. The spin structure on $W$ extends to a spin
structure on $W^\#$. We choose the boundary spin structure on $M^\#$.
If $p$ and $q$ are in different connected
components of $M$, then the spin structure on $W^\#$ (resp. $M^\#$) 
is uniquely determined.
However, if $p$ and $q$ are in the same connected component, then
there are two non-equivalent choices of spin structures on $W^\#$ 
that extend the spin structure on $W$. These two spin structures 
arise from each other by a diffeomorphism as explained in the following lemma.

\begin{lemma}\label{spin.struc.unique} Let $n=\dim M =\dim W -1\geq 2$.
If $\kappa_1$ and $\kappa_2$ are two spin structures on $W^\#$ such that 
the inclusions
$(W,\kappa_j) \hookrightarrow (W^\#,\kappa_j^\#)$, 
$j=1,2$, are spin preserving embeddings, then there is a diffeomorphism
$f:W^\#\to W^\#$ which is the identity 
on~$W$ and such that $f^*\kappa_1=\kappa_2$.
Similarly, if $\chi^\#_1$ and $\chi^\#_2$ are two spin structures on $M^\#$ 
such that the inclusions 
$(M\setminus (B_p(\ep)\cup B_q(\ep)),\chi)\hookrightarrow 
(M^\#,\chi_j^\#)$, $j=1,2$, are spin preserving embeddings, 
then there is a diffeomorphism
$f:M^\#\to M^\#$ which is the identity 
on~$M\setminus (B_p(\ep)\cup B_q(\ep))$ and such that $f^*\chi_1=\chi_2$.
\end{lemma}

The proof of this lemma is straightforward.

Also note a particularity in the case $n=2$:
$S^1\times [-1,1]$ carries $2$ non-equivalent spin structures.
However, it is clear from the construction above that
the restriction of $\chi^\#$ to
$S^1\times [-1,1]$ is the product spin structure of the bounding
spin structure on $S^1$ with the unique spin structure on $[-1,1]$.

\subsection{Identifying spinors for different metrics} \label{identify}
Throughout the paper, we need to identify spinor fields for a manifold $M$
with a fixed (topological) spin structure $\chi$, 
but two different metrics $g$ and $h$.
In order to recall this identification we follow 
\cite{bourguignon.gauduchon:92}.
 
Let $x \in M$. Since the metrics $g$ and $h$ are   
symmetric and positive definite, there is a unique symmetric 
$b^g_h  \in \End(T_x M)$ such that  for all $v,w \in T_x M$, 
$$h(b^g_h v ,b^g_h w) = g(v, w).$$ Note that $b^g_h$ depends smoothly on $x$.
Hence, the map 
\begin{eqnarray*}
  (b^g_h)^n:{\mathrm{SO}(M,g)}&\to&{\mathrm{SO}(M,h)}\\
    (e_1, \cdots,e_n)&\mapsto&(b^g_h(e_1), \cdots ,b^g_h(e_n))
\end{eqnarray*}
is an isomorphism of $SO(n)$-principal bundles
between the oriented frame bundles of $(M,g)$ and $(M,h)$. 
This map commutes with the right action of $\mathrm{SO}_n$.
Furthermore, $(b^g_h)^n$ can be lifted to    
\dgARROWLENGTH=2em   
$$\begin{diagram}   
\node{\mathrm{Spin}(M,g)}\arrow[2]{e,t}{\beta^g_h}\arrow{s}\node[2]{\mathrm{Spin}(M,h)}\arrow{s}\\   
\node{\mathrm{SO}(M,g)}\arrow[2]{e,t}{(b^g_h)^n}\arrow{s}\node[2]{\mathrm{SO}(M,h)}\arrow{s}\\   
\node{M}\arrow[2]{e,t}{\mathrm{Id}}\node[2]{M}   
\end{diagram}$$   
Hence, we obtain a map, also denoted by $\beta^g_h$, 
between the spinor bundles $\Sigma_g M$ and $\Sigma_h M$ 
in the following way:   
\begin{eqnarray}\label{indentspinbun}   
\Sigma_g M=\mathrm{Spin}(M,g)\times_\sigma\Sigma_n 
&\longrightarrow& \Sigma_h M=\mathrm{Spin}(M,h)\times_\sigma\Sigma_n\nonumber\\   
\psi=[s,\varphi]&\longmapsto& \beta^g_h{\psi} =[\beta^g_h(s),\varphi]   
\end{eqnarray} 
where $(\sigma,\Sigma_n)$ is the complex spinor representation, and 
$[s,\varphi]$ denotes the equivalence class of $(s,\varphi)$ under the 
diagonal action of $\Spin(n)$.
The identification $\beta^g_h$ of spinors preserves 
the pointwise norm of spinors. Apparently,
$\beta^g_h\circ \beta^h_g=\id$. However, in general, for three scalar products
$g$, $k$, $h$ we have $\beta^h_g\circ \beta^k_h \circ \beta^g_k\neq \id$.
It is a direct consequence of the construction of the bundle
$\mathrm{Spin}(M,g)$ and the spinor bundle $\Si_g M$ that 
there is a section $B^g_h:T^*M\otimes\Sigma_g M\to \Sigma_g M$
\begin{eqnarray}\label{connect.eqn}
\beta^g_h{\nabla^g_X \psi} - \nabla^{h}_X \beta^g_h{\psi}
  = \beta^g_h(B^g_h(X \otimes\psi))
\end{eqnarray} 
The expression $B^g_h$ only depends on $g$, $h$ and their first derivatives.
In particular, $\|B^g_h\|\to 0$ if $h$ converges to $g$ 
in the $C^1$-topology.  
The Dirac operators with respect to $g$ and $h$ 
are locally defined as 
$D_g = \sum_i e_i \cdot \nabla_{e_i}^g$ and 
$D_h= \sum_i b^g_h(e_i) \cdot \nabla_ {b^g_h( e_i)}^h$. It follows that
\begin{eqnarray} \label{comparison_D}
|\beta^g_h{D_g \psi} - D_{h} \beta^g_h{\psi}| \leq C^g_h 
\left(|\psi|+ |\nabla^g \psi| \right)
\end{eqnarray} 
 where $C^g_h\to 0$ if $h\to g$ in $C^1$.

In the particular case where the metric $h$ is conformal to $g$, i.e.\ when
we have $h=f^2g$ where $f$ is a smooth positive function, then by
\cite{hitchin:74,hijazi:86} one has that 
\begin{eqnarray} \label{conf_D}
D_h( f^{-\frac{n-1}{2}} \beta^g_h\psi) = f^{-\frac{n+1}{2}} \beta^g_h{D_g \psi}.\end{eqnarray}

This equation implies, in particular, that the dimension of the kernel 
of the Dirac operator is constant on any conformal class.

\subsection{The first eigenvalue of the Dirac operator in a conformal
  class} \label{variational}
Let $(M,g,\chi)$ be a compact spin manifold of dimension $n$, $\ker D=\{0\}$.
For  $\psi\in\Gamma(\Sigma M)$ we define    
$$J(\psi)=\frac{\Big(\int_M|D\psi|^{\frac{2n}{n+1}}\,dv_g\Big)^\frac{n+1}{n}}{\int_M <D\psi,\psi>\,dv_g}.$$  
Using techniques from \cite{lott:86},
Ammann proved in \cite{ammann:03} that    
\begin{eqnarray} \label{funct}   
\lamin(M,g,\chi)=\inf_\psi J(\psi)   
\end{eqnarray}   
where the infimum is taken over the set of smooth spinor fields for which    
$$\left(\int_M  <D\psi,\psi>\,dv_g \right)>0.$$   

We will need the following result: 
\begin{prop}
Let $(M,g,\chi)$ be a compact spin manifolds of dimension $n \geq 2$.      
Then, 
\begin{equation} \label{aubin}   
  \lamin(M,g,\chi) \leq \lamin(\mS^n) = \frac{n}{2}\, \om_n^{\frac{1}{n}} 
\end{equation} 
where $\om_n$ stands for the volume of the standard sphere   
$\mS^n$.
\end{prop}

The proposition was proven in \cite{ammann:03} using geometric
methods if $n\geq 3$. In the case $n=2$ the article \cite{ammann:03}
only provides a proof if $\ker D=\{0\}$.
Another method, that yields the proposition in full generality,
is to construct for any $p\in M$ and $\ep>0$ 
a suitable test spinor field $\psi_\ep$ supported in $B_p(\ep)$ 
that verifies 
  $$J(\psi_\ep)\leq \lamin(\mS^n)+o(\ep)$$
(see \cite{ammann.humbert.morel:p03a,grhu} for details).

If inequality~\eref{aubin} holds even strictly, 
then one can show that the infimum in equation~\eref{funct} is attained, say
in $\phi$. Then the infimum in the definition of $\lamin(M,g,\chi)$ 
is attained in the generalized conformal metric
$\ti g:=|\phi|^{4/(n-1)}g$ (see \cite{ammann:comin} for details).
 
This fact, summarized in the following 
theorem will be a central ingredient in the proof of our main result 
Theorem~\ref{main}.

\begin{theorem}[\cite{ammann:comin,ammann.habil}] \label{attain}
Assume that Inequality (\ref{aubin}) is strict. Then there exists a 
spinor field $\phi \in C^{2,\alpha}(\Si M) \cap   
C^{\infty}(\Si M \setminus \phi^{-1}(0))$ on $(M,g)$, 
$\al\in(0,1)\cap (0,2/(n-1)]$ such that 
\begin{equation}\label{eq.dir.con}   
  D(\phi)= \lamin(M,g,\chi)\; |\phi|^{{2 \over n-1}} \phi, \hbox{ and }   
{\parallel \phi \parallel}_{{2n \over n-1}}=1
\end{equation} 
Furthermore, there is a generalized conformal metric $\tilde{g}$ (see
\cite{ammann.habil} for the definition)  such   
that $\la_1^+(\tilde{g}) \Vol(M,\tilde{g})^{{1/ n}}= \lamin(M,g,\chi)$. 
\end{theorem}

\subsection{A conformal Sobolev inequality}

Note that we have the pointwise
  $$|D\psi|\leq \sqrt{n}|\na \psi|$$
which implies
  $$\frac{\int_{M}  |D \psi|^{\frac{2n}{n+1}} \,dv_g}
 {\int_{M} |\nabla \psi|^{\frac{2n}{n+1}}}\leq \sqrt{n}$$
Elliptic regularity theory states that this fraction is also bounded from below.
\begin{prop} \label{sobolev}
Let $(M,g,\chi)$ be a compact spin manifold of dimension $n \geq 2$. We assume
that the Dirac operator $D$ is invertible. Then, there
exists a constant $C >0$ such that for all  spinor fields $\psi$ of class
$C^1$, we have :
\begin{eqnarray} \label{sob1} 
\int_{M}  |\na \psi|^{\frac{2n}{n+1}} \,dv_g 
\leq C
 \int_{M}  |D \psi|^{\frac{2n}{n+1}} \,dv_g 
\end{eqnarray}
and 
\begin{eqnarray} \label{sob2mod} 
{ \left( \int_{M}  |\psi|^{\frac{2n}{n-1}} \,dv_g \right)}^\frac{n-1}{n}
\leq C
{ \left( \int_{M}  |D \psi|^{\frac{2n}{n+1}} \,dv_g \right)}^\frac{n+1}{n}.
\end{eqnarray}
\end{prop}

\proof{} Inequality (\ref{sob1}) is a classical elliptic inequality for
invertible elliptic operators on compact manifolds. 
It is equivalent to the fact that $D^{-1}$
is a continuous operator from $L^{\frac{2n}{n+1}}$ to $W^{1,\frac{2n}{n+1}}$ 
(see e.g.\ \cite{taylor:81}). Inequality (\ref{sob2mod}) 
is the classical Sobolev embedding theorem which asserts that the Sobolev
space $W^{1,\frac{2n}{n+1}}(M)$ is continuously embedded into
$L^{\frac{2n}{n-1}}(M)$.
 
\begin{corollary}
Let $(M,g,\chi)$ be a compact Riemannian spin manifold with invertible 
Dirac operator. 
Let $(\Omega,h)$ be conformal to an open subsets of $(M,g)$.
Then there is a constant $C=C(M,g,\chi)$ such that any compactly
supported spinor field $\psi\in\Gamma(\Sigma(\Omega,h,\chi))$ of class $C^1$ 
satisfies 
$${ \left( \int_{\Omega}  |\psi|^{\frac{2n}{n-1}} \,dv_h \right)}^\frac{n-1}{n}
\leq C
{ \left( \int_{\Omega}  |D_h \psi|^{\frac{2n}{n+1}} \,dv_h
  \right)}^\frac{n+1}{n}.$$
\end{corollary}

\proof{}
Let $f$ be the conformal factor, i.e. $h=f^2g$. 
With the notation of paragraph  \ref{identify}, we let 
$\phi:= f^{\frac{n-1}{2}} \beta^h_g\psi \in
\Gamma(\Sigma(\Omega,g,\chi))$.
Since the map $\beta^h_g:\Si_hM\to\Si_gM$ is
a pointwise isometry, and since $dv_h= f^n\, dv_g$ and by
(\ref{conf_D}), we get that 
    $$\int_{\Omega}  |\psi|^{\frac{2n}{n-1}} \,dv_h = \int_{\Omega}
|\phi|^{\frac{2n}{n-1}} \,dv_g$$
and 
$$ \int_{\Omega}  |D_h \psi|^{\frac{2n}{n+1}} \,dv_h  = \int_{\Omega}  |D_g \phi|^{\frac{2n}{n+1}} \,dv_g.$$  
The result is then an immediate consequence of inequalities (\ref{sob1}) and 
(\ref{sob2mod}).

\qed
 
\begin{example} \label{cyl}
Let $M=\mS^n$ be the sphere. 
Then the Mercator projection is a conformal embedding 
$F:\mR\times \mS^{n-1}\to \mS^n$ 
$$F :\begin{pmatrix}t\cr x_1 \cr \cdots \cr x_n \end{pmatrix}\mapsto 
  \begin{pmatrix}\tanh(t)\cr
                  {x_1 \over \cosh( t) }\cr 
                  \cdots \cr  
                  {x_n \over\cosh(t) }
  \end{pmatrix}.$$
The image of $F$ is $\mS^n$ with the North and South pole removed.
Recall from the preliminaries that $\mS^n$ is always equipped 
with the spin structure $\chi^n$ that arises as restriction
of the unique spin structure on $\overline{B_{0,\mR^n}(1)}$ to the boundary.
In the case $n\geq 2$ this is the unique spin structure, in the case 
$n=1$ there are two possible spin structures on $S^1$, one of them is the 
bounding spin structure $\chi^1$.
One easily verifies that $F^*(\chi^n)$ is the product structure
of the unique spin structure on $\mR$ and of $\chi^{n-1}$.

Hence, any compactly supported spinor $\psi$ on the cylinder
$\mR\times \mS^{n-1}$ satisfies
\begin{equation}\label{cyl.estim}
{ \left( \int_\mR\int_{\mS^{n-1}}  |\psi|^{\frac{2n}{n-1}} \,dv_{\ground{n-1}} 
\,dt\right)}^\frac{n-1}{n}
\leq C
{ \left( \int_\mR\int_{\mS^{n-1}}  |D \psi|^{\frac{2n}{n+1}} \,dv_{\ground{n-1}}\,dt \right)}^\frac{n+1}{n}.
\end{equation}
\end{example}

\section{Approximation by locally flat manifolds}\label{sec.approx}

The first step in the proof of Theorem~\ref{main} is to approximate $(M,g)$ 
by metrics that are flat in a neighborhood of $p$ and $q$, and then show that 
$\lamin$ does not change much under this approximation.
The results of this section will allow us to assume 
that $g$ is flat in a neighborhood of $p$ and $q$.

\begin{lemma} \label{app_flat}
Let $(M,g, \chi)$ be a compact spin manifold of dimension $n \geq 2$. We
assume that the Dirac operator $D$ is invertible.  Let   $p,q \in M$. 
There is a family of Riemannian metrics 
$(g_\ep)_{\ep\in(0,\al)}$ on $M$, such that each $g_\ep$ 
flat in a neighborhood of $p$ and $q$ and such that $g_\ep \to g$  
in $C^1$ when $\ep \to 0$.  
\end{lemma}

\proof{}
Let $\eta:[0,\infty)\to [0,1]$ be a smooth function that equals to $1$ on 
$[0,1]$ and whose support is contained in $[0,2)$. 
Choose a small number $\ep>0$.
On $M\setminus (B_p(2\ep)\cup B_q(2\ep))$ we define $g_\ep=g$.
In normal coordinates centered in $p$ and defined on $B_p(2\ep)$ we write  
$g(x)=\sum_{ij}g_{ij}(x)\,dx^i\,dx^j$. We define on $B_p(2\ep)$
  $$g_\ep(x):=\sum_{ij}\Bigl(\eta(\ep^{-1}|x|)\,\de_{ij} + (1-\eta(\ep^{-1}|x|))\,g_{ij}(x)\Bigr) \,dx^i\,dx^j,$$ 
and similarly on $B_q(2\ep)$. Then one calculates that
$g_\ep \to g$  in $C^1$ when $\ep \to 0$. 
\qed 

\begin{proposition} \label{conv}
Let $(M,g,\chi)$ be a compact spin manifold of dimension $n \geq 2$
with invertible Dirac operator $D$.
If $g_{\ep}$ is a sequence of Riemannian metrics converging to $g$ in
$C^1$, then
$$\lim_{\ep\to 0}\lamin(M,[g_{\ep}],\chi)= \lamin(M,[g],\chi).$$
\end{proposition}

{\bf Proof of Proposition~\ref{conv}.\\}
Since $D$ is invertible, and as the spectrum of $D$ seen as 
of function of $g$ depends continuously on $g$, it follows  
that the Dirac operator on $(M,g_\ep,\chi)$, denoted by $D_{\ep}$, 
is invertible as well.  Let $J$ (resp. $J_{\ep}$) be the 
functional (see Paragraph \ref{variational}) associated to
$(M,[g],\chi)$ (resp. $(M,[g_{\ep}],\chi)$).
Let $\psi \in \Gamma (\Sigma_g M)$ be a smooth spinor field on $(M,g,\chi)$.
We denote by $D_\ep$ the
Dirac operator acting on $(M,g_\ep,\chi)$.We can choose $\psi$ such that $J(\psi) \leq \lamin(M,[g],\chi) + \de$
where $\delta>0$ is small. By 
\eref{comparison_D}, \eref{sob1} and\eref{sob2mod}
one easily gets
that $$\limsup_{\ep \to 0} \lamin(M,[g_{\ep}],\chi) \leq 
 \limsup_{\ep \to 0} J_\ep(\beta^g_{g_\ep} \psi) = J(\psi) \leq
\lamin(M,[g],\chi) + \de.$$ 
Since $\delta>0$ is arbitrary, we obtain that 
\begin{eqnarray} \label{lim+}
\lim_{\ep \to 0} \lamin(M,[g_{\ep}],\chi) \leq \lamin(M,[g],\chi).
\end{eqnarray}
Now, let $\psi_\ep \in \Gamma (\Sigma_g M)$ be a smooth spinor field on
$(M,g,\chi)$ such that 
$$J_\ep(\beta^g_{g_\ep} \psi_\ep ) \leq \lamin(M,[g_{\ep}],\chi) + \de.$$
In order to abbreviate, we set $\overline\psi_\ep:=\beta^g_{g_\ep} 
\psi_\ep$.
Without loss of generality, we can assume that 
$$\int_M \<D_\ep \overline{\psi}_\ep, \overline{\psi}_\ep\> \,dv_{g_\ep} =1.$$
Then, since $J_\ep(\overline{\psi}_\ep )$ is bounded (by
Relation (\ref{aubin})), then 
$\int_M |D_{\ep} \overline{\psi}_\ep|^\frac{2n}{n+1} \,dv_{g_{\ep}}$ 
is also bounded.
We have
$$\left| \left( 
\int_M |D_{\ep} \overline{\psi}_\ep|^\frac{2n}{n+1} \,dv_{g_{\ep}}  \right)^\frac{n+1}{2n}  - 
 \left(  \int_M  |D \psi_\ep|^\frac{2n}{n+1} \,dv_g  \right)^\frac{n+1}{2n}
\right| 
\leq  
\left| \left( 
\int_M |D_{\ep} \overline{\psi}_\ep|^\frac{2n}{n+1} \,dv_{g_\ep}  \right)^\frac{n+1}{2n}  - 
 \left(  \int_M  |D_{\ep} \overline{\psi}_\ep|^\frac{2n}{n+1} \,dv_g  \right)^\frac{n+1}{2n}
\right| $$
\begin{eqnarray} \label{i1}
+  \left| \left( 
\int_M |D_{\ep} \overline{\psi}_\ep|^\frac{2n}{n+1} \,dv_g  \right)^\frac{n+1}{2n}  - 
 \left(  \int_M  |D \psi_\ep|^\frac{2n}{n+1} \,dv_g  \right)^\frac{n+1}{2n}
\right|.
\end{eqnarray}
By definition of $g_{\ep}$, it is clear that 
$$ \lim_{\ep \to 0} \left| \left( 
\int_M |D_{\ep} \overline{\psi}_\ep|^\frac{2n}{n+1} \,dv_{g_\ep}  \right)^\frac{n+1}{2n}  - 
 \left(  \int_M  |D_{\ep} \overline{\psi}_\ep|^\frac{2n}{n+1} \,dv_g  \right)^\frac{n+1}{2n}
\right|= 0.  $$
In addition, we get from relation (\ref{comparison_D}) that 
 
$$\left| \left( 
\int_M |D_{\ep} \overline{\psi}_\ep|^\frac{2n}{n+1} \,dv_g  \right)^\frac{n+1}{2n}  - 
 \left(  \int_M  |D \psi_\ep|^\frac{2n}{n+1} \,dv_g  \right)^\frac{n+1}{2n}
\right| \leq a_{\ep}  \left(  \int ( |\psi_\ep| + |\nabla_g \psi_\ep|)
  |^\frac{2n}{n+1} dv_g \, \right)^\frac{n+1}{2n}.$$
where  $\lim_{\ep \to 0} a_{\ep} = 0$.
Using the triangle inequality and H\"older's inequality, 
we see the existence of a constant $C >0$
independent of $\ep$ such that 
 $$\left| \left( 
 \int_M |D_{\ep} \overline{\psi}_\ep|^\frac{2n}{n+1} \,dv_g  \right)^\frac{n+1}{2n}  - 
  \left(  \int_M  |D \psi_\ep|^\frac{2n}{n+1} \,dv_g  \right)^\frac{n+1}{2n}
 \right| \leq C  a_{\ep} \left(  \left(  \int  |\psi_\ep|^\frac{2n}{n-1}
     dv_g \, \right)^\frac{n-1}{2n}+   \left(  \int |\nabla_g
     \psi_\ep|^\frac{2n}{n+1} dv_g \, \right)^\frac{n+1}{2n} \right).$$
 By inequalities (\ref{sob1}) and (\ref{sob2mod}), we obtain that 
 $$\left| \left( 
 \int_M |D_{\ep} \overline{\psi}_\ep|^\frac{2n}{n+1} \,dv_g  \right)^\frac{n+1}{2n}  - 
  \left(  \int_M  |D \psi_\ep|^\frac{2n}{n+1} \,dv_g  \right)^\frac{n+1}{2n}
 \right| \leq C  a_{\ep}  \left(  \int |\nabla_g \psi_\ep|^\frac{2n}{n+1}
   dv_g \, \right)^\frac{n+1}{2n}.$$
 Coming back to (\ref{i1}), this clearly implies that 
 $\int_M  |D \psi_\ep|^\frac{2n}{n+1} \,dv_g $ is bounded and that $\int_M |\psi_\ep|^{\frac{2n}{n-1}}\, \,dv_g $ is bounded, too.
Since $\lim_{\ep \to 0} a_\ep = 0$, this also implies that 
\begin{eqnarray} \label{limnum} 
\lim_{\ep \to 0} \left| 
\int_M |D_{\ep} \overline{\psi}_\ep|^\frac{2n}{n+1} \,dv_{g_{\ep}} - 
\int_M  |D \psi_\ep|^\frac{2n}{n+1} \,dv_g \right|  = 0.
\end{eqnarray}
Similarly, using relation (\ref{comparison_D}) we have 
\begin{eqnarray} \label{limden}
\lim_{\ep \to 0} \left|   \int_M \<D_\ep \overline{\psi}_\ep , \overline{\psi}_\ep \>
  \,dv_{g_\ep } -  \int_M \<D \psi_{\ep}, \psi_\ep\> \,dv_g \right| = 0.
\end{eqnarray}
Relations (\ref{limnum}) and (\ref{limden}) imply that 
 $$\lamin(M,[g_{\ep}],\chi) + \delta  \geq 
 J_\ep(\overline{\psi}_\ep) \geq   J(\psi_\ep) + o(1)  \geq
\lamin(M,[g],\chi)+o(1).$$ 
Since $\de$ is arbitrary, we get that 
\begin{eqnarray} \label{lim-}
\liminf_{\ep \to 0} \lamin(M,[g_{\ep}],\chi) \geq \lamin(M,[g],\chi).
\end{eqnarray}
Together with (\ref{lim+}), this proves Theorem \ref{conv}.

\section{Construction of the metrics on $M^\#$} \label{metric}

The aim of this section is to construct the sequence of metrics $g_\ep^\#$
of Theorem~\ref{main}. 
Using the result of the previous section, 
\emph{we can assume from now on and in the rest of the article 
that the metric $g$ is flat on $B_p(\rho)$ and $B_q(\rho)$ for a small 
$\rho>0$.} 
For $0<\al<\be<\rho$ we introduce the notation
  $$ B_{p,q}(\al) = B_p(\al) \cup B_q (\al)  \hbox{ and }
   C_{p,q}(\al,\be) = B_{p,q}(\be)\setminus B_{p,q}(\al).$$
Let $\ep >0$ be small. We explained in the introduction that $M^\#$ is obtained
(as a topological space) by gluing a cylinder $[-1,1]\times S^{n-1}$ with
$M \setminus B_{p,q}(\ep)$ along $\partial  B_p(\ep)$ on one side and 
along $\partial  B_q(\ep)$ on the other side. Evidently this can also be 
expressed by saying  
  $$M^\#= M \setminus B_{p,q}(\ep)/\sim$$
where $\sim$ indicates that we glue $\pa B_p(\ep)$ with $\pa B_q(\ep)$ via an
orientation preserving diffeomorphism. $M^\#$ is equipped with a differential
structure and a spin structure such that 
$I_\ep:M \setminus B_{p,q}(\ep)\to M^\#$
is smooth and compatible with the spin structures. 
We also introduce for all $a\in (\ep,\rho)$ the notation 
  $$H(a,\ep):=I_\ep(C_{p,q}(\ep,a))=C_p(\ep,a)\cup C_q(\ep,a)/\sim.$$
Hence $M^\#$ is the disjoint union of $M\setminus B_{p,q}(a)$ and 
$H(a,\ep)$. 

Let us denote by $d_p$ (resp. $d_q$) the distance in $(M,g)$ 
to the point $p$ (resp. $q$).
We define a function $f_{\ep}$ on 
$(M \setminus B_{p,q}(3\ep)) \cup C_{p,q}(\ep,2\ep)$
by setting 
\[ f_\ep(x)  =  \left| \begin{array}{ccc}
1 & \hbox{ if } & x \in 
M \setminus
B_{p,q}(3 \ep) \\
d_p^{-1} & \hbox{ if } & x \in C_p(\ep,2 \ep) \\
d_q^{-1} & \hbox{ if } & x \in C_q(\ep,2 \ep). \\
\end{array} \right. \]
We can extend $f_\ep$  smoothly and positively 
to $M\setminus B_{p,q}(\ep)$ such that $f_\ep$ 
satisfies on $C_{p,q}(2\ep,3\ep)$:
$$|\na f_{\ep}| \leq \frac{2}{\ep}.$$
As $(C_p(\ep,2\ep),f_\ep^2 g)$ and $(C_q(\ep,2\ep),f_\ep^2 g)$ are
isometric to $([0,\log 2)\times S^{n-1},dt^2+\ground{n-1})$,
$f_\ep^2g$ defines a metric $g_\ep^\#$ on $M^\#$, or more precisely:
there is a unique metric $g_\ep^\#$ on $M^\#$ such that
  $$f_\ep^2g =I_\ep^* g_\ep^\#\qquad \mbox{holds on }M\setminus B_{p,q}(\ep).$$

Note that $(H(a,\ep),g_\ep^\#)$ is then
conformal (but in general not isometric) 
to $((-\log(a/\ep),\log(a/\ep))\times S^{n-1},dt^2+\ground{n-1})$.

\section{Proof of Theorem \ref{main}}\label{sec.proof}

In this section we will prove that 
the metrics $(g_{\ep}^\#)_\ep$ are the desired metrics.
We denote by $D_{\ep}$ the Dirac operator acting on 
$(M^\#,g_{\ep}^\#,\chi^\#)$. We set $\la= \lamin(M,g,\chi)$ and 
$\la_{\ep} = \lamin(M^\#,g_{\ep}^\#,\chi^\#)$.
We denote by $J$ (resp. $J_{\ep}$)  the functional
associated to $\la$ (resp. $\la_{\ep}$) (see Paragraph \ref{variational}).

\begin{lemma}
\begin{eqnarray}\label{g>ge}
\limsup_{\ep \to 0} \la_{\ep} \leq \la.
\end{eqnarray}
\end{lemma}

\proof{}
Let $\gamma >0$ be small and let  $\psi$ be a smooth spinor field such that
$J(\psi) \leq  \la + \gamma$. Clearly, for each small number $\al >0$, 
one can construct a cut-off function $\eta_{\al} \in
C^{\infty}(M)$, such that $0 \leq \eta_{\al} \leq 1$,  
equal to $1$ on $M \setminus B_{p,q}(2 \al)$, equal to $0$ on 
$B_{p,q}(\al)$ and which satisfies $|\nabla \eta_{\al}| \leq
\frac{2}{\al}$. As easily seen, we have $\lim_{\al \to 0}J(\eta_{\al}
  \psi) = J(\psi)$.
We choose $\al$ small enough such that 
$J(\eta_{\al}
  \psi) \leq \la + 2 \gamma$. Now, if $3 \ep < \al$, then $\eta_{\al}
  \psi$ is supported on the common part of $(M,g,\chi)$ and $(M^\#,
  g_{\ep}^\#,\chi^\#)$ and hence can be seen as a spinor field 
on $(M^\#, g_{\ep}^\#,\chi^\#)$. 
We  have 
  $$\la_{\ep}\leq  J_{\ep}(\eta_{\al}  \psi)  \leq \la +2 \gamma$$
for all small $\ep>0$.
Since $\gamma$ is arbitrary, we proved the lemma.\qed

It remains to prove that 
\begin{eqnarray} \label{g<ge}
\liminf_{\ep \to 0} \la_{\ep} \geq \la.
\end{eqnarray}
This inequality is more involved than Equation\eref{g>ge} and will occupy 
the rest of this section.

We set $\la_0:=\liminf_{\ep \to 0}\la_{\ep}$ and pass to a sequence of 
$\ep_i\to 0$ with $\lim_{i\to \infty} \la_{\ep_i}=\la_0$.
To simplify notation we write $\la_i:=\la_{\ep_i}$, 
$D_i:=D_{\ep_i}$, $g_i^\#:=g_{\ep_i}^\#$,
$dv_i:= dv_{g_{\ep_i}^\#}$, and so on. In the following arguments
we will frequently pass to subsequences. 
Slightly abusing the notation, we will continue with the same 
index notation $\la_i$, $D_i$, $g_i^\#$, and so on.

In the case $\la_0 = \lamin(\mS^n)$
Equation~\eref{g<ge} follows directly from \eref{aubin}.
Hence, after possibly passing to a subsequence, we can assume that 
$\la_i < \lamin(\mS^n)$ for all $i\in \mN$. 
By Theorem~\ref{attain}, there exists a spinor
field $\psi_i$ of class $C^1$ defined on  
$(M^\#, g_i^\#,\chi^\#)$ that satisfies
\begin{eqnarray} \label{equation}
D_i \psi_i = \la_i |\psi_i|^{\frac{2}{n-1}} \psi_i
\end{eqnarray}
with the normalization 
\begin{eqnarray} \label{normal}
\int_{M^\#}  |\psi_i|^{\frac{2n}{n-1}} \,dv_i = 1.
\end{eqnarray}

\begin{lemma}
For any $\delta >0$ and $a_0>0$ there is an $a\in (0,a_0)$ such that after 
possibly passing to a subsequence 
\begin{eqnarray} \label{liminf}
\liminf_{i\to \infty }
\int_{C_{p,q}(a,2a)} |\psi_i|^{\frac{2n}{n-1}} \,dv_g
\leq  \de
\end{eqnarray}
and such that $g$ is flat on $C_{p,q}(a,2a)$. 
\end{lemma}

Note that the integral above has a meaning  for $\ep_i \leq a$, since
$C_{p,q}(a,2a ) \subset M^\#$. If in addition, $3 \ep_i <a$, we have 
$g_i=g$ on $C_{p,q}(a,2a)$.
For the proof of  (\ref{liminf}) we proceed by contradiction. 
We assume that
(\ref{liminf}) is false. Let $N_0$ be such that  $g$ is flat on $B_{p,q}(2^{-N_0})$. For each $N$, we can find a sequence $\ep_i\to 0$
such that for all $k \in \{N_0, \cdots, N\}$, and
we have 
$$ \lim_{i \to \infty}
\int_{C_{p,q}( 2^{-(k+1)},  2^{-k} )} |\psi_i|^{\frac{2n}{n-1}} \,dv_g
>  \de.$$
Since 
$$\int_{M^\#}  |\psi_i|^{\frac{2n}{n-1}} \,dv_i=1$$
and since the $C_{p,q}( 2^{-(k+1)},  2^{-k} )$ are disjoint, we obtain a
contradiction if $N - N_0 > \frac{1}{\de}$ for $a= 2^{-(N+1)}$.
This proves Relation (\ref{liminf}). \\

\noindent We fix $N$ for which  Relation (\ref{liminf}) is verified.
Eventually extracting a subsequence of $(\psi_i)$, we can assume that
the limit infimum in (\ref{liminf}) is a limit. 
In other words, we have found a number $a$ such that
\begin{eqnarray} \label{liminf2}
\int_{C_{p,q}(a,2a)} |\psi_i|^{\frac{2n}{n-1}} \,dv_g \leq \delta.
\end{eqnarray}

Since 
$$\int_{M \setminus B_{p,q}(2a)}  |\psi_i|^{\frac{2n}{n-1}} \,dv_g + 
\int_{C_{p,q}(a,2a)} |\psi_i|^{\frac{2n}{n-1}} \,dv_g +
\int_{H(a,\ep)}
|\psi_i|^{\frac{2n}{n-1}} \,dv_i = 1$$
since for $\ep$ small 
$$ \int_{C_{p,q}(a,2a)} |\psi_i|^{\frac{2n}{n-1}} \,dv_g \leq \de$$
and since $\de$ can be chosen lower than $\frac{1}{3}$, we have for large $i$:
\begin{eqnarray} \label{concM}
\limsup_{i\to \infty}\int_{M \setminus B_{p,q}(2a)}  |\psi_i|^{\frac{2n}{n-1}} \,dv_i \geq \frac13
\end{eqnarray}
or
\begin{eqnarray} \label{concMM}
\limsup_{i\to \infty}\int_{H(a,\ep_i)}
|\psi_i|^{\frac{2n}{n-1}} \,dv_i\geq \frac13.
\end{eqnarray}
where $H(a,\ep_i)$ is defined as in section \ref{metric}.

The theorem now follows from the following two lemmata.

\begin{lemma}\label{lem.conc.M}
Relation (\ref{concM}) implies $\la_0\geq \la$.
\end{lemma}

\begin{lemma}\label{lem.conc.handle}
Relation (\ref{concMM}) implies $\la_0\geq \lamin(\mS^n)$.
\end{lemma}

\proof{of Lemma~\ref{lem.conc.handle}}\\
We suppose that Inequality~\eref{concM} holds.

{\bf Step 1. $\{i\in\mN\,|\, \la_i= 0\}$ is finite.}

We prove this step by contradiction and assume that 
after passing to a subsequence $\la_i = 0$ for all $i$.

This means that the spinors $\psi_i$ are harmonic on $(M^\#,g_i^\#)$.
 
By Example \ref{cyl}, the manifold $(H(2a,\ep), g_i^{\#})$
is conformal to a subdomain of the sphere $\mS^n$. 
By Proposition~\ref{sobolev} there exists a constant $C$ 
which does not depend on $i$ 
and $\delta$  such
that for all spinor field $\theta$ of class $C^1$ and whose support is
included  in $H(2a,\ep_i)$, we have 
\begin{eqnarray}\label{sob12}
{ \left( \int_{M^\#}|\theta|^{\frac{2n}{n-1}} \,dv_i
  \right)}^\frac{n-1}{n} \leq  C
{ \left( \int_{M^\#}  |D_i \theta|^{\frac{2n}{n+1}} \,dv_i \right)}^\frac{n+1}{n}.
\end{eqnarray}
Let $\eta \in C^{\infty}(M^\#)$,
$0 \leq \eta \leq 1$  be a cut-off function
equal  to $1$ on $H(a,\ep_i)$ supported in $H(2a,\ep_i)$,
and which satisfies $|\nabla \eta | \leq \frac{2}{a}$.
Then we apply Inequality \eref{sob12} with $\theta = \eta \psi_i$.
Since $\psi_i$ is harmonic and by \eref{concMM}, we get that 
$$3^{-\frac{n-1}{n}}  \leq C{ \left( \int_{M^\#} |(\nabla \eta) \cdot  \psi_i|^{\frac{2n}{n+1}} \,dv_i \right)}^\frac{n+1}{n}.$$
Using the H\"older inequality, we get that 
$$3^{-\frac{n-1}{n}} \leq  C
{ \left( \int_{M^\# } |\nabla \eta|^n  \,dv_i \right)}^\frac{2}{n}
{ \left( \int_{Supp(\nabla \eta)}  |\psi_i|^{\frac{2n}{n-1}} \,dv_i
\right)}^\frac{n-1}{n}.$$
Since $|\nabla \eta| \leq \frac{2}{a}$, since $\nabla \eta$ is supported
in $C_{p,q}(a,2a)$ whose volume for $g$ is bounded by $C a^n$ and since $g_i$ is constant equal to $g$ on $C_{p,q}(a,2a)$, 
there exists a constant $C'$ independent of $i$ and
$\delta$ such that 
\begin{eqnarray} \label{nablaeta}
\int_{M}  |\nabla \eta|^n  \,dv_i \leq C'.
\end{eqnarray}
Using (\ref{liminf2}), we get that 
$$\int_{Supp(\nabla \eta)}  |\psi_i|^{\frac{2n}{n-1}} \,dv_i 
\leq \delta.$$
We obtain
$$\frac{1}{3^\frac{n-1}{n} } \leq  C{ C'}^\frac{2}{n} \delta^\frac{n-1}{n}.$$
If $\delta$ is small enough, we get a contradiction. This proves Relation 
step 1.\\

{\bf Step 2. Proof of the inequality
\begin{eqnarray} \label{limla}
\lim_{i\to \infty} \la_i \geq  \lamin(\mS^n).
\end{eqnarray}}\\
Let $\eta$ be the same function as above.
Since the manifold $(H(2a,\ep), g_i^{\#})$
is conformal to a subdomain of the sphere $\mS^n$ and since $\eta \psi_i$ 
is supported in $H(2a,\ep_i)$, we have 
\begin{eqnarray} \label{quot}
\lamin(\mS^n) \leq J(\eta \psi_i).
\end{eqnarray}
By
Equation~\eref{equation}, we have 
$$\left( \int_{M^\#}\<D_i (\eta \psi_i) , \eta \psi_i\> \,dv_i \right)
= \la_i \int_{M^\#} \eta^2 |\psi_i|^\frac{2n}{n-1}\,dv_i \, \in \mR.$$
In consequence, since  
$$ \underbrace{\int_{ M^\#}\<D_i (\eta \psi_i),\eta \psi_i\> \,dv_i}_{\in \mR} = 
\underbrace{\int_{M^\#} \<\nabla \eta \cdot \psi_i, \eta \psi_i\> \,dv_i}_{\in i\mR} + 
\underbrace{\int_{M^\#} \eta^2 \<D_i \psi_i, \psi_i\> \,dv_i}_{\in \mR},$$
we get that $\int_{M^\#}\<\nabla \eta \cdot \psi_i, \eta
\psi_i\>\,dv_i=0$. Since $\eta \equiv 1$ on $H(a,\ep_i)$, we obtain that 
\begin{eqnarray} \label{denom}
 \left( \int_ {M^\#}\<D_i (\eta \psi_i) , \eta \psi_i\> \,dv_i \right)
\geq \la_i \int_{H(a,\ep_i)} |\psi_i|^\frac{2n}{n-1} \,dv_i.
\end{eqnarray}
We also have by Equation (\ref{equation}) 
$${ \left( \int_ {M^\#}|D (\eta \psi_i)|^\frac{2n}{n+1} \,dv_i \right)
}^\frac{n+1}{n}  =  { \left( \int_ {M^\#}\left|\nabla \eta \cdot \psi_i + 
\eta\,\la_i|\psi_i|^\frac{2}{n-1} \psi_i \right|^\frac{2n}{n+1} \,dv_i\right)
}^\frac{n+1}{n}.$$
Again since $\eta \equiv 1$ on $H(a,\ep_i)$, we can write that 
$$\int_{M^\#}\left|\nabla \eta \cdot \psi_i + \eta\,\la_i 
|\psi_i|^\frac{2}{n-1} \psi_i \right|^\frac{2n}{n+1} \,dv_i  \leq $$
$$\int_{C_{p,q}(a,2a)} \left|\nabla \eta \cdot \psi_i + \eta\,\la_i 
|\psi_i|^\frac{2}{n-1} \psi_i \right|^\frac{2n}{n+1} \,dv_i + 
\la_i^\frac{2n }{n+1} \int_{H(a,\ep_i)}|\psi_i|^\frac{2n }{n-1}\,dv_i.$$
Using the fact that for $s,t \in \mR$, we have 
$|s+t|^\frac{2n}{n+1} \leq |2s|^\frac{2n}{n+1}+ |2t|^\frac{2n}{n+1}$ and
using H\"older inequality, we
get that 
$$ 
\int_{C_{p,q}(a,2a)} 
\left|\nabla \eta \cdot \psi_i + \eta\,\la_i 
|\psi_i|^\frac{2}{n-1} \psi_i \right|^\frac{2n}{n+1} \,dv_i
 $$
$$\leq 2^\frac{2n}{n+1} \int_{C_{p,q}(a,2a)}  |\nabla \eta|^\frac{2n}{n+1}
|\psi_i|^\frac{2n}{n+1} \,dv_i +
2^\frac{2n}{n+1}  \la_i^\frac{2n}{n+1} 
\int_{C_{p,q}(a,2a)}|\psi_i|^\frac{2n}{n-1} \,dv_i$$
$$ \leq  2^\frac{2n}{n+1} \left(  \int_{C_{p,q}(a,2a)}  |\nabla \eta|^n  \,dv_i\right)^\frac{2}{n+1} \left(  \int_{C_{p,q}(a,2a)}  |\psi_i|^\frac{2n}{n-1} \,dv_i  
  \right)^\frac{n+1}{n-1} + 2^\frac{2n}{n+1} \,\la_i^\frac{2n}{n+1}
  \delta. $$
Using again relation (\ref{nablaeta}), we get that 

$$\int_{C_{p,q}(a,2a)} 
\left|\nabla \eta \cdot \psi_i + \eta\,\la_i 
|\psi_i|^\frac{2}{n-1} \psi_i \right|^\frac{2n}{n+1} \,dv_i \leq  C (\delta^\frac{n+1}{n-1} + \delta)$$
where $C$ is a constant independent of $i$ and $\delta$.
Finally, we obtain that 
\begin{eqnarray} \label{num}
\left( \int_ {M^\#}|D \eta \psi_i|^\frac{2n}{n+1} \,dv_i \right)^\frac{n+1}{n} 
\leq \left( C (\delta^{\frac{n+1}{n-1}} + \delta ) 
     + \la_i^{\frac{2n}{n+1}}  \int_{H(a,\ep_i)} 
       |\psi_i|^{\frac{2n}{n-1}}\,dv_i\right)^{\frac{n+1}{n}}.
\end{eqnarray}
Plugging (\ref{num}) and (\ref{denom}) into (\ref{quot}), we get that 
$$\lamin(\mS^n) \leq \frac{{ \left( C (\delta^\frac{n+1}{n-1} + \delta ) + 
\la_i^{\frac{2n}{n+1}}  \int_{H(a,\ep_i)} 
|\psi_i|^\frac{2n}{n-1}\,dv_i \right) }^\frac{n+1}{n}}
{\la_i \int_{H(a,\ep_i)} |\psi_i|^\frac{2n}{n-1} \,dv_i}.$$
Since $$\int_{H(a,\ep_i)} |\psi_i|^\frac{2n}{n-1} \,dv_i \geq \frac{1}{3}$$
and letting $\delta$ go to $0$, we obtain 
Relation (\ref{limla}).
Since by (\ref{aubin}), we have $\la \leq \lamin(\mS^n)$, we get 
\eref{g<ge}.

This proves the lemma.

\proof{of Lemma~\ref{lem.conc.M}}
We just sketch the proof of  Lemma~\ref{lem.conc.M} because the method 
is the same as the proof of Lemma~\ref{lem.conc.handle},
At first, we prove step 1. 
We proceed by contradiction and assume that $\la_i = 0$ for infinitely many
$i\in \mN$. 
This means that the spinors $\psi_i$ are harmonic on $(M^\#,g_i^\#)$. 
By Proposition~\ref{sobolev}, there exists $C>0$ such that for all spinor 
$\theta$ of class $C^1$ defined on $M$, we have 
$${ \left( \int_{M}  |\theta|^{\frac{2n}{n-1}} \,dv_{g_i^\#} \right)}^\frac{n-1}{n}
\leq 
C { \left( \int_{M}  |D  \theta|^{\frac{2n}{n+1}} \,dv_g  \right)}^\frac{n+1}{n}.$$
We apply this inequality with $\theta= (1- \eta)  \psi_i$ where 
$\eta$ is defined as in the proof of Lemma~\ref{lem.conc.handle}
and we obtain a contradiction if
$\delta$ is small enough. This proves step 1.
Then, we say that 
$$\la \leq J((1- \eta) \psi_i).$$
As in the proof of Lemma~\ref{lem.conc.handle}, 
we obtain (\ref{g<ge}) in the limit $\de\to 0$. 
This proves Theorem~\ref{main}.

\section{The Dirac operator on Riemann surfaces}\label{riem.surf}

In this section we want to prove Theorem~\ref{main.twodim}.
In fact, we use Theorem~\ref{main} to calculate the $\tau$-invariant for all 
compact Riemann surfaces equipped with spin structures. 

As already mentioned in the introduction, the Atiyah-Milnor-Singer 
invariant $\alpha$ is a ring homomorphisms
from the ring of spin-cobordism classes into $KO^{-*}(pt)$. It allows to
define an index theorem for the Dirac operator 
that is non-trivial in dimensions
$0,1,2,4\ \mod 8$ (see e.g.~\cite[II.7]{lawson.michelsohn:89}).

Let us to recall an equivalent definition of
the $\alpha$-invariant in the special case $n=2$. In this case 
$\alpha(M)\in KO^{-2}\{pt\}=\mZ_2:=\mZ/(2\mZ)=\{0,1\}$. 
Some more details can also be found in \cite{kusner.schmitt:p96} and
\cite[Section~2 and 3]{ammann.baer:02}. 
We will also recall the index theorem in dimension $2$, 
Theorem~\ref{theo.ind.two}.

A quadratic form is a map
$q:H_1(M,\mZ_2)\to \mZ_2$ such that   
$$q(a+b)= q(a)+q(b) + a\cap b$$
holds for all $a,b\in H_1(M,\mZ_2)$. Here $\cap$ 
denotes the intersection form 
$H_1(M,\mZ_2)\times  H_1(M,\mZ_2)\to \mZ_2$ which is a non-degenerate
(anti-)symmetric bilinear form on $H_1(M,\mZ_2)$.

The difference of two such quadratic forms is a linear map
$H_1(M,\mZ_2)\to \mZ_2$, and vice versa if one adds a linear
map $H_1(M,\mZ_2)\to \mZ_2$ to a quadratic form, one easily sees
that one obtains a quadratic form again. 
The space of quadratic forms on $H_1(M,\mZ_2)$ is an
affine spaces modelled 
on the space $H^1(M,\mZ_2):=\Hom(H_1(M,\mZ_2),\mZ_2)=\Hom(H_1(M,\mZ),\mZ_2)$.

We will now associate to any spin structure $\chi$ on a Riemann surface~$M$ 
a quadratic form
$q_\chi:H_1(M,\mZ_2)\to \mZ_2$. This association will define
a bijection from the set of equivalence classes of spin structures
to the set of quadratic forms.

For simplicity of notation, 
we fix a Riemannian metric $g$ on $M$. Let $\SO(M,g)$
denote the $S^1$-principal bundle of positively oriented orthonormal frames
on $M$. If one specializes the description of 
a spin structure (Subsection~\ref{subsec.spin}) to dimension $2$, then
a spin structure $\chi$ consists of a principal $S^1$-bundle 
$\Spin(M,g,\chi)\to M$ and a double covering
$\mu_\chi:\Spin(M,g)\to \SO(M,g)$ with 
$\mu_\chi(\phi\cdot z)= \mu_\chi(\phi)\cdot z^2$ for all 
$\phi\in \Spin(M,g)$ and all $z\in S^1$.

Any homology class in $H^1(M,\mZ_2)$ can be represented by a closed embedded 
loop $\gamma:\mR/(L\mZ)\to M$, $L>0$, parameterized by arclength.
If $\dot \gamma:\mR/(L\mZ) \to \SO(M,g)$ lifts to a 
map $\mR/(L\mZ)\to \Spin(M,g)$,
then we define $q_\chi(\gamma):=1$, otherwise we define $q_\chi(\gamma):=0$.
One checks that if $\gamma_1$ and $\gamma_2$ represent the same homology class
in $H_1(M,\mZ_2)$, then $q_\chi(\gamma_1)=q_\chi(\gamma_2)$, 
hence $q_\chi$ defines a map $q_\chi:H_1(M,\mZ_2)\to \mZ_2$. 
One checks that this map is in fact a quadratic map.

For any quadratic map $q:V\to \mZ_2$ associated to a non-degenerate
symmetric bilinear form on a finite-dimensional 
$\mZ_2$-vector space $V$ one defines 
one defines the Arf-invariant
  $$\Arf(q):=\frac{1}{\sqrt{\# V}}\,
    \sum_{a\in V}(-1)^{q_\chi(a)}.$$
One can check that the sum is either $+1$ or $-1$. 
We now define $\al(M,\chi)$ via 
$(-1)^{\al(M,\chi)}=\Arf(q_\chi)$.  

\begin{example}
Let $M=S^2$ with the spin structure $\chi^2$. 
Then $q_{\chi^2}:\{0\}\to\mZ_2$, $q_{\chi^2}(0)=0$, $\al(S^2,\chi^2)=0$.
\end{example}
\begin{example}
Let $M$ be of genus $1$, i.e.\ diffeomorphic to $T^2$, with a metric $g$. 
After performing
a conformal change (which does not affect neither the quadratic form nor the 
dimension of the kernel) we can assume that $g$ is flat, i.e.\ $(M,g)$ is
isometric to $\mR^2/\Gamma$ equipped with the euclidean metric
for a lattice $\Gamma\subset \mR^2$, $\Gamma=\pi_1(M)=H_1(M,\mZ)$. 
Then $SO(M,g)$ is a trivial $S^1$ bundle, where a trivialization is given by a
parallel frame $(e_1,e_2)$. We can also write 
$SO(M,g)=(\mR^2/\Gamma)\times \SO(2)$.
Any group homomorphism 
$\gamma:\Gamma\to \{-1,+1\}\subset \ker(\Spin(2)\to\SO(2))\subset
\Spin(2)=S^1$ 
defines a diagonal action of 
$\Gamma$ on $\mR^2\times \Spin(2)$, and we obtain a 
$\Spin(2)$-principal bundle by factoring out this action 
 $$\Spin(M,g)_\Gamma:=\mR^2\times_\ga \Spin(2).$$
This principal bundle together the natural map 
\begin{eqnarray*}
\chi_\ga:\Spin(M,g)=\mR^2\times_\ga \Spin(2)&\to& \SO(M,g)=(\mR^2/\Gamma)\times \SO(2)\\
  {}[(x,z)]_\ga &\mapsto & (x+\Gamma,z^2)
\end{eqnarray*}
defines a spin structure $\chi_\ga$ on $M$. Note that in the sense of
$\Spin(2)$-principal bundles, the possible $\Spin(M,g)$ are all
equivalent.
However, $\chi_{\ga_1}$ and 
$\chi_{\ga_2}$ are equivalent in the sense of spin structures
iff $\ga_1=\ga_2$. Furthermore any spin structure
on $M$ is obtained in this way.  Denote the image of $v\in \Gamma$ in 
$H_1(M,\mZ_2)=\Gamma\otimes_\mZ \mZ_2$ by $\ol v$. 
The quadratic form $q$ of the spin structure associated to $\gamma$
then fulfills
  $$\begin{matrix}
    q(\ol v)=\gamma(v)+1\hfill& \mbox{for $\ol v\neq 0$,}\hfill\\
    q(\ol v)=\gamma(v)=0\hfill& \mbox{for $\ol v= 0$.}\hfill    
 \end{matrix}
$$ 
This implies $\al(M,\chi)=1$ if $\chi$ is the spin structure associated to the 
trivial map $\gamma$, and $\al(M,\chi)=0$ in all other cases.
The boundary of a solid torus has a map
$\gamma$ which is non-trivial, hence $\al(M,\chi)=0$ in this case.

On the other hand we can also calculate the dimension of the kernel of $D$.
As $g$ is flat, $D\phi=0$ is equivalent to $\na\phi=0$.

If $\gamma$ is the trivial map, then the spinor bundle is trivialized
by parallel spinors, i.e.\ $\dim \ker D_g=2$. If $\gamma$ is non-trivial,
then $\dim \ker D_g=0$. 

This terminates the example, and we return to the general case.
\end{example}

We know recall the index theorem for Dirac operators on 
compact Riemann surfaces. Note that
the spinor bundle $\Sigma M\to M$ carries the structure of a quaternionic
vector bundle over $M$, and the quaternionic multiplication commutes with 
the Dirac operator. Hence, the complex dimension of any eigenspace
of the Dirac operator is divisible by $2$.

\begin{theorem}[Index theorem]\label{theo.ind.two}
For any compact surface with 
Riemannian metric $g$ and spin structure $\chi$ we have
  $${(\dim_\mC \ker D) /2}\equiv \al(M,\chi)\mod 2.$$ 
\end{theorem}
In the examples above, we have verified this relation  
if $M$ is diffeomorphic to $S^2$ or $T^2$.
We will sketch a short proof of the theorem in Remark~\ref{ind.proof}.

\begin{remark} According to \cite{hitchin:74} the (complex) 
dimension $k_g$ of the 
kernel of the Dirac operator on a compact Riemann surface $(M,g)$
of genus $\gamma$ is at most $\gamma+1$. 
Hence, $k_g$ is already
determined by $\al(M,\chi)$ and $\gamma$ if $\gamma\leq 2$, or if 
$\gamma=3$, $\alpha(M,\chi)=1$. However, in all other cases, $k_g$ depends
on the conformal class of $g$ (see
\cite{hitchin:74}, \cite{baer.schmutz:92}). The spectrum of 
the Dirac operator depends continuously on $g$ in the
$C^1$-topology. Hence, $k_g\geq \limsup_{h\to g} k_h$.
If $g_i\to g$ with $k_{g_i}<k_g$ 
then due to the symmetric of the spectrum of $D$,
some positive and some negative eigenvalue converges to $0$, both 
having the same, even multiplicity. 
Hence, if $k_g$ jumps then by a multiple
of $4$, and $k_g\mod 4$ is therefore independent of $g$.
\end{remark}

\proof{of Theorem~\ref{main.twodim}}
In the case $\al(M,\chi)=0$, the index theorem implies that the 
Dirac operator has a kernel for any given metric on $M$. 
This immediately implies $\tau(M,\chi)=0$ which yields
the first case in Theorem~\ref{main.twodim}.

In order to derive the second case, we have to study the behavior
of the Arf-invariant on products and the $\al$-invariant under connected sum.

If 
$q_i:V_i\to \mZ_2$ are quadratic maps for $i=1,2$, then on
$V_1\oplus V_2$ we have a product quadratic map defined via
$q_1\otimes q_2:V_1\oplus V_2\to \mZ_2$, 
$(q_1\otimes q_2)(v_1,v_2)=q_1(v_1)+q_2(v_2)$. One checks that
$\Arf(q_1\otimes q_2)=\Arf(q_1) \Arf(q_2)$.

The quadratic form of a disjoint union of $(M_1,\chi)$ and $(M_2,\chi_2)$ 
is just the product quadratic map $q_{\chi_1}\otimes q_{\chi_2}$. 
Furthermore, one easily sees that the quadratic forms 
of $(M_1\# M_2,\chi_1\#\chi_2)$ and $(M_1,\chi_1)\dot\cup (M_2,\chi_2)$ 
can be identified.  
 
It follows that 
  $$\al(M_1\# M_2,\chi_1\#\chi_2)\equiv \al((M_1,\chi_1)\dot\cup
  (M_2,\chi_2))\equiv \al(M_1,\chi_1)+\al(M_2,\chi_2)\mod 2$$
As $\alpha(S^2,\chi_2)=0$, it follows that 
$\al(M,\chi)$ is a spin-cobordism invariant. Hence, $\al$ defines a map 
from the $2$-dimensional spin-cobordism group to $\mZ_2$. 

Inversely, it can be shown (e.g. with the statements of 
\cite[section~3]{ammann.baer:02}):

\begin{lemma}
Let $M$ be an compact oriented surface of genus $\gamma$ 
with spin structure $\chi$, and $\al(M,\chi)=0$.
Then $(M,\chi)$ is spin-diffeomorphic to the connected sum of $\gamma$ 
$2$-tori 
  $$T^2\# T^2 \# \cdots\# T^2$$
where each $2$-torus carries a spin structure that is associated
to a non-trivial homomorphism. 
\end{lemma}
 
\begin{corollary}
Let $M$ be an compact oriented surface of genus $\gamma$ 
with spin structure $\chi$, and $\al(M,\chi)=0$.
Then $(M,\chi)$ is obtained from $(S^2,\chi^2)$ by a sequence of 
$0$-dimensional surgeries.
\end{corollary}

Hence, Theorem~\ref{main} implies that 
$\tau(M,\chi)\geq \tau(S^2,\chi^2)=2\sqrt{\pi}$,
but as $\tau(M,\chi)\leq 2\sqrt{\pi}$ (Proposition~\ref{aubin}), we obtain
the conclusion in the second case of Theorem~\ref{main.twodim}.\qed

\begin{remark}\label{ind.proof}
With the  methods provided in this section 
Theorem~\ref{theo.ind.two} can be proved in a simple and geometric way.
From the construction of $\al(M,\chi)$ out of the quadratic form
it is clear that $\al(M,\chi)$ is preserved under $0$-dimensional
surgery, in particular
it is additive under connected sum. Together with $\alpha(S^2,\chi^2)=0$
it follows that it is  a spin-cobordism invariant.
We have already seen that $\dim \ker D\mod 4$ is invariant on the metric.
The
following Proposition implies that $\dim \ker D$ might only jump by multiples
of $4$ when a $0$-dimensional surgery is performed that 
introduces a (long and thin) cylinder. Hence, $\dim \ker D\mod 4$ is
spin-cobordant as well. As we have already checked the index theorem for
tori, the index theorem Theorem~\ref{theo.ind.two} follows in general. 

\begin{proposition}
Let $(M,g,\chi)$ be a compact spin manifold with Dirac operator $D$, 
and let $(M^\#,\chi^\#)$
be obtained from $M$ by $0$-dimensional surgery. Then there is a sequence
of metrics $g_i^\#$, $i\to \infty$, on $M^\#$ such that the Dirac
operator $D_i$ on $(M^\#,g_i^\#,\chi^\#)$ satisfies
\begin{enumerate}[{\rm (a)}]
\item $\dim \ker D_i$ is independent on $i$,  
\item $\dim \ker D_i\leq \dim \ker D$,
\item If $\dim \ker D_i<\dim \ker D$, then there are positive and 
negative eigenvalues converging to $0$. Their combined multiplicity is equal
to  $\dim \ker D-\dim \ker D_i$. 
\end{enumerate}   
\end{proposition}
Part (a) and (b) of the Proposition are proved in 
\cite{ammann.dahl.humbert:p06}. Part (c) is not proved in there
explicitly, but the
arguments and constructions in \cite{ammann.dahl.humbert:p06} 
can be adapted such that 
we obtain (c).

The metrics $g_i$ are unchanged outside the attached cylinders, 
and the cylinders
equipped with the metrics $g_i$ get longer and thinner when 
$i$ tends to $\infty$.
\end{remark}


\vspace{1cm}               
Authors' address:               
\nopagebreak   
\vspace{5mm}\\   
\parskip0ex     
\vtop{   
\hsize=8cm\noindent   
\obeylines               
Bernd Ammann and Emmanuel Humbert               
Institut \'Elie Cartan BP 239              
Universit\'e de Nancy 1               
54506 Vandoeuvre-l\`es -Nancy Cedex               
France                           
\vspace{0.5cm}               
               
E-Mail:               
{\tt bernd.ammann at gmx.net and humbert at iecn.u-nancy.fr}               
}   

                    
\end{document}